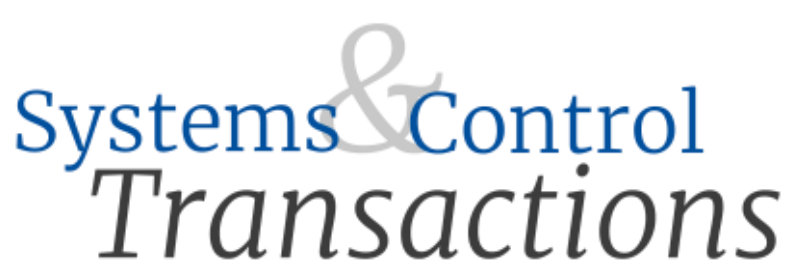

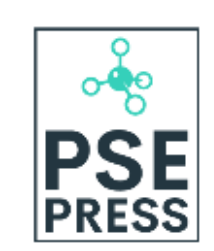

# REVISITING SIMULTANEOUS METHODS FOR DYNAMIC OPTIMIZATION IN THE GPU ERA

**Joseph W. Choi[a] and Sungho Shin[a]***

[a] Massachusetts Institute of Technology, Department of Chemical Engineering, Cambridge, Massachusetts, USA
* Corresponding Author: sushin@mit.edu

## ABSTRACT

We revisit the classical topic in dynamic optimization: sequential vs simultaneous methods for solving DAE-constrained optimization problems, with a particular focus on how graphics processing unit (GPU) computing changes their effectiveness. Sequential methods offer key advantages through adaptive time stepping at the differential–algebraic equation (DAE) solver level, which is especially effective for handling stiff systems. However, long-time-horizon simulations remain a computational bottleneck, as time integration is inherently sequential and limits parallelization within the optimization algorithm. In contrast, simultaneous approaches are well-suited for parallel computing. They address these limitations by exploiting the highly repetitive structure of discretized DAE systems at the function evaluation level and leveraging sparse linear algebra routines that enable elimination tree–level parallelism. Although simultaneous methods typically lack adaptive time stepping, this limitation can often be mitigated by choosing a sufficiently fine initial mesh or iteratively adjusting mesh coarseness in an outer loop. In this work, we revisit the capabilities of the simultaneous approach in a GPU computing environment and assess its performance against a sequential method baseline. We employ a simultaneous approach based on orthogonal collocation within an open-source modeling framework and apply it to parameter estimation benchmarks from systems biology. We evaluate both GPU and CPU solvers on the simultaneous formulation. Our results show that, although less reliable, the simultaneous approach achieves up to $5.4\times$ speedup compared to the sequential baseline among the largest instances where both methods solve successfully. The advantage of the simultaneous method with respect to problem size is more pronounced on GPUs than on CPUs.

## INTRODUCTION

The emergence of graphics processing units (GPUs) acceleration techniques in nonlinear programming (NLP) warrants revisiting the choice between sequential and simultaneous methods for solving differential-algebraic equation (DAE)-constrained optimization problems. In the sequential method, the DAE solution profile is obtained by a *sequential* forward solve at each iterate [1]. While this sequential time stepping enables the solver to overcome stiff dynamics, it also prevents parallelization within the optimization algorithm [1]. In contrast, the simultaneous method discretizes the DAE as algebraic constraints, enabling the state profiles and remaining free variables to be solved and optimized *simultaneously* as a single, large-scale NLP [1].

The cost of solving NLPs is dominated by two computational steps: sparse linear algebra to solve the Newton step at each iteration, and automatic differentiation (AD) to supply the first- and second-order derivatives. While this cost typically scales with NLP size, recent work has mitigated some of these issues by parallelizing certain aspects of both steps on the GPU. In particular, MadNLP.jl, a filter line-search interior point solver, enables sparse linear algebra routines on the GPU by reformulating the KKT system into a pivoting-free factorizable form [2]. The majority of the benefit from GPU parallelization is obtained at this sparse linear algebra level [2]. ExaModels.jl, an algebraic modeling framework, parallelizes AD through a *single instruction, multiple data* (SIMD) abstraction [2]. Under this abstraction, the NLP is specified as a few expression patterns, each a "single instruction", with the AD performed in parallel over an array of "multiple data". Combined, GPU-accelerated sparse linear algebra and SIMD abstraction are especially beneficial for exploiting NLPs with sparse and highly repetitive structures.

Although the resulting NLP from the simultaneous method is large, most of its size stems from the fine collocation mesh required to accurately capture the governing dynamics. The constraints themselves however are sparse and consist of only a few distinct expression patterns. As a result, the simultaneous method is particularly well suited for SIMD abstraction. Furthermore, challenges associated with mesh fineness are typically handled by computationally demanding adaptive mesh refinement or error estimation strategies [1]. SIMD abstraction may potentially remediate mesh refinement issues by supporting the adoption of an overly generous initial mesh at little computational expense, closing the gap on the advantages of the sequential method. This context motivates revisiting the simultaneous method to answer the question: does GPU acceleration improve the efficacy of the simultaneous method, and if so, under what circumstances is it most effective?

To compare the simultaneous and sequential methods, we consider parameter estimation problems from systems biology because they cover a broad range of transient behaviors and problem cases. The Benchmarking Initiative's PEtab collection [3] gathers such problems in a standardized format referred to as *PEtab models* [4]. This collection pairs with PEtab.jl [5], a sequential method implementation for solving PEtab models in Julia. Because ExaModels.jl and MadNLP.jl are also based in Julia, a GPU-based simultaneous approach can be compared within the same environment. However, while PEtab.jl provides the baseline for the sequential method, a GPU-compatible simultaneous method implementation for solving PEtab models remains to be built.

### Related works

Several works implement the automatic formulation of NLPs based on the simultaneous method, but the following are the most relevant. InfiniteExaModels.jl [6], a backend for InfiniteOpt.jl [7], implements the simultaneous method in a GPU-compatible fashion by automatically transcribing infinite-dimensional problems in ExaModels. KIPET [8] applies the simultaneous method to solve parameter estimation problems from chemical reaction kinetics, but is implemented in Pyomo using Pyomo.DAE [9], which lacks native support for GPU-based NLP solvers.

### Contributions

This work makes two contributions. First, we introduce ExaModelsPEtab.jl, a GPU-compatible simultaneous method analog to PEtab.jl that automatically formulates an NLP in ExaModels.jl directly from PEtab models. Our work differs from InfiniteExaModels.jl by specializing the construction of GPU kernels for the structure of PEtab models, and from KIPET by enabling the solution entirely on the GPU. Second, we conduct a benchmark comparing the sequential PEtab.jl and simultaneous ExaModelsPEtab.jl formulations across GPU and CPU solvers. To our knowledge, no such comparison exists between the sequential and simultaneous methods in the context of GPU computing.

### Overview

This paper is structured as follows: **Sequential and Simultaneous Methods for Parameter Estimation** introduces the sequential and simultaneous methods applied to parameter estimation problems in a general form. **Implementation of ExaModelsPEtab** details the core implementation features of ExaModelsPEtab.jl in the context of SIMD abstraction. **Numerical Results** presents the benchmark results, highlighting the respective advantages and limitations of both methods. **Conclusion** provides concluding remarks and areas for future directions.

## SEQUENTIAL AND SIMULTANEOUS METHODS FOR PARAMETER ESTIMATION

This section reviews the sequential and simultaneous methods for solving parameter estimation problems. First, we introduce the estimation problem in a general form. We then present the resulting reduced and large-scale NLP formulations arising from the sequential and simultaneous methods, respectively.

### Parameter estimation problem formulation

For simplicity, we consider ordinary differential equation (ODE)-constrained estimation problems of the general form

$$\min_{z(\cdot),\, \theta^{\mathrm{L}} \leq \theta \leq \theta^{\mathrm{U}}} \Phi(z(\cdot),\, \theta) \tag{1a}$$

$$\text{s.t.}\ \frac{dz^c}{dt} = f^c(z^c(t), \theta, t), \quad t \in T^c \tag{1b}$$

$$z^c(t=0) = z_0^c(\theta) \tag{1c}$$

$$c = 1, \dots, N_c \tag{1d}$$

where $c$ is the index for different experiment runs, $N_c$ is the number of experiments, $z^c(\cdot)$ is the state profile, $\theta$ is the estimated parameter vector with bounds $\theta^{\mathrm{L}}$ and $\theta^{\mathrm{U}}$, $\Phi(\cdot)$ is a cost function that penalizes the model and measurement deviation, $f^c(\cdot)$ is the right-hand side function, $t \in T^c$ denotes time over an interval, and $z_0^c(\cdot)$ is the initial condition function. The superscript $c$ denotes a condition-specific profile or mapping. Hence, the objec-tive depends on $z(\cdot)$, the state profile across all times and conditions. We omit algebraic constraints for simplicity, noting that observable expressions are defined implicitly within $\Phi(\cdot)$.

## Sequential method formulation

This section briefly describes the sequential method for solving Equation (1) to provide context for a later defined benchmark metric. Instead of optimizing over the states, an embedded ODE solver integrates eq. (1b) at each iteration of the optimization cycle, reducing the NLP to

$$\min_{\theta^{\mathrm{L}} \leq \theta \leq \theta^{\mathrm{U}}} \Phi_{\text{seq}}(\theta)$$

where $\Phi_{\text{seq}}(\cdot)$ is the sequential method cost function which only explicitly depends on $\theta$. Since the integration is performed at the ODE solver level, an accurate state profile can be obtained to a desired tolerance. Step directions, which are only required for $\theta$, are then computed using first- and second-derivative information from either direct or adjoint sensitivity methods [1].

## Simultaneous method formulation

To solve Equation (1) with the simultaneous method, the states are discretized on a finite collocation mesh and kept as decision variables, resulting in a large-scale NLP. The collocation mesh is defined by the collocation points $t_{ij} = t_{(i-1)0} + h_i \tau_j$ where $h_i$ is the width of interval $i$ for $i = 1, \ldots, N$ intervals, and $\tau_j \in [0,1], j = 0, \ldots, K$ are interpolation points for a total $K+1$ points per interval. Given a collocation mesh, the NLP is given by

$$\min_{z,\, \theta^{\mathrm{L}} \leq \theta \leq \theta^{\mathrm{U}}} \Phi_{\text{sim}}(z,\, \theta) \tag{2a}$$

$$\text{s.t.} \sum_{j=0}^{K} a_j\,(\tau_k) z_{ij}^c = h_i f^c(z_{ik}^c, \theta, t_{ik}),$$

$$i = 1, \ldots, N, k = 1, \ldots, K \tag{2b}$$

$$z_{(i+1)0}^c = \sum_{j=0}^{K} b_j\, z_{ij}^c, \quad i = 1, \ldots, N-1 \tag{2c}$$

$$z_{10}^c = z_0^c(\theta) \tag{2d}$$

$$c = 1, \ldots, N_c \tag{2e}$$

where $z_{ij}^c$ are the discretized states, $z$ is all of the discretized states, $\Phi_{\text{sim}}(\cdot)$ is the simultaneous method cost function which depends explicitly on both $z$ and $\theta$, and $a_j(\tau_k)$, $b_j$ are weights that depend on the chosen basis polynomials. In contrast to the sequential method, step directions are computed for $z$ and $\theta$ since they are solved and optimized simultaneously.

# IMPLEMENTATION OF EXAMODELSPETAB

This section details the key feature implementations of ExaModelsPEtab.jl in the context of SIMD abstraction. In particular, we present our design choices regarding the collocation constraint formulation and provide a code snippet to illustrate its construction in ExaModels.jl.

## Collocation equation with condition independence

We adopt the orthogonal collocation formulation with $N$ intervals, degree $K$ Lagrange basis polynomials, and Gauss-Legendre roots. We also remove the experimental-condition dependence from the right-hand side function by lifting all condition-dependent variables to the function-input space. Applied to eq. (2b), the collocation constraints become

$$\sum_{j=0}^{K} \frac{d\ell_j^K(\tau_k)}{d\tau} z_{ij}^c = h_i f(z_{ik}^c, \theta, t_{ik};\, p^c),$$
$$i = 1, .., N,\ k = 1, \ldots, K,\ c = 1, \ldots, N_c \tag{3}$$

where $p^c$ are condition-dependent variables, $f(\cdot\,;\, p^c)$ is the lifted condition-independent right-hand side function, and the weights $d\ell_j^K / d\tau$ are the $\tau$ derivative of the $j$th Lagrange basis polynomial of degree $K$.

By writing the collocation equations in the form Equation (3), the number of unique expression patterns is reduced to the number of state variables $N_z$ from $N_c \times N_z$. The cost of this lifting is at most an additional $\dim p^c$ often extremely simple expression patterns. This resulting algebraic structure aligns well with the core design philosophy of SIMD abstraction: few expression patterns are repeated across all experimental conditions and collocation points, presenting the opportunity for massive parallel AD evaluations. The same principle is applied to eq. (2a), eq. (2c), and eq. (2d) to obtain similar parallelization.

## Constructing GPU kernels with ExaModels.jl

The following listing demonstrates how the GPU AD kernels for Equation (3) are constructed in ExaModels.jl:

```
coll_cons = ExaModels.@add_con(core,
   -hi*f[v](z[:,i,k,c],theta[:],tij,p[:,c])
   for (i, k, c, hi, tij) in iter_coll)
ExaModels.@add_con!(core, coll_cons,
   (i,k,c) => dljdtauk*z[v,i,j,c]
   for (i, j, k, c, dljdtauk) in iter_coll!)
```

Here, `@add_con` creates the constraint set for a single right-hand side function $f_v(\cdot)$, while `@add_con!` augments each constraint with the summed terms over every interpolation point. The parallel AD evaluations are enabled by the "multiple data" iterators `iter_coll` and `iter_coll!` which contain all of the indices, functional inputs, and constant terms.

# NUMERICAL RESULTS

This section details the benchmark setup and results comparing both GPU and CPU solvers with the simultaneous method using ExaModelsPEtab.jl against the sequential method using PEtab.jl. The main purpose of this comparison is to examine the discussed GPU opportunities for exploiting the structure of the simultaneous method, rather than to declare one method superior over the other. The code for reproducing our results can be found in https://github.com/mit-shin-group/exa-petab-focapo-2027.

## Benchmark setup

### Hardware

The benchmarks were performed on a single machine with two 16-core Intel Xeon Gold 6130 CPUs, 1 TB of RAM, and two NVIDIA Quadro GV100 GPUs with 32 GB each.

### Selected benchmark models

Only 20 out of the 35 total PEtab models from the Benchmarking Initiative's PEtab collection [3] were chosen for benchmarking. The remaining 15 models were excluded because they either contain the model feature "Possible Discontinuities", which is not yet supported by ExaModelsPEtab.jl, or failed to compile using PEtab.jl, which our package relies on. Still, the 20 models cover a broad range of biological dynamics: 11 intracellular signaling and regulatory networks, 2 metabolic and enzyme kinetics, 3 infection and immune models, and 4 population-level epidemiological models [3]. Problem sizes span from 3 to 75 ODEs and 3 to 113 estimated parameters which span up to 14 orders of magnitude.

### Optimization framework configurations

**ExaModels.jl + MadNLP.jl (GPU/CPU):** ExaModels.jl v0.11.1 with MadNLP.jl v0.10.1. Both the GPU and CPU use the LiftedKKT formulation [10]. Only the linear solver differs between the two configurations; GPU uses cuDSS (CUDSS.jl v0.7.0) with MadNLPGPU v0.10.1 and CUDA.jl v6.2.0, and CPU uses ma57. Both optimizers use $10^{-6}$ tolerance for relative KKT error and $10^{-4}$ acceptable tolerance within 15 iterations. We use degree $K = 4$ polynomials because our experiments showed it offered the most consistency between performance and collocation accuracy. We obtain both a fixed mesh and initial guesses for the discretized states by solving the ODE system at the PEtab model-provided nominal parameter values $\theta_0$ using the same ODE solver and configurations as PEtab.jl, which we describe below.

**PEtab.jl + Best optimizer (CPU):** PEtab.jl v5.3.0 with Optim.jl v2.1.0, Fides.jl v1.2.1, and OrdinaryDiffEq.jl v7.0.0. Defining equivalent optimizers or configurations between the sequential and simultaneous methods is impossible due to their fundamentally different formulations. To make a comparison nonetheless, we report the best performance across three solvers recommended by PEtab.jl for small, medium, and large models respectively [5]: `Optim.IPNewton` with an exact `ForwardDiff` hessian, `Fides.CustomHessian` with `GaussNewton`, and `Fides.BFGS`. For each, we adopt the default optimizer configurations as chosen by PEtab.jl. Notably, the default gradient tolerance is $10^{-6}$, which we used to base our KKT tolerance choice for MadNLP. At the ODE solver level, we adopt the integrator and its configurations defaulted by PEtab.jl's model-size heuristic. All integrators default to $10^{-8}$ absolute and relative tolerance. Exact integrator configurations based on this heuristic can be found within the PEtab.jl default options documentation [5].

**Both configurations:** All solution frameworks share uncapped iteration count, one hour solve and model compile caps, and PEtab model-provided nominal estimated parameter values as initial guesses. Any unspecified configuration implies default settings, which can be found on the respective package repositories.

### Performance metrics

We evaluate two performance metrics across the solution configurations. The shifted geometric mean solve time $t_{\text{SGM}} = \left(\prod_{i=1}^{n}\left(t_{\text{solve},i} + 0.01\right)\right)^{1/n} - 0.01$ for $n = 5$ re-solves is reported in seconds, where $t_{\text{solve},i}$ is the optimizer wall time for run $i$. We compare solution quality with the relative objective gap (ROG) given by

$$\text{ROG} := \frac{\left|\Phi_{\text{seq,PEtab}}(\theta^*_{\text{Exa}}) - \Phi_{\text{seq,PEtab}}(\theta^*_{\text{PEtab}})\right|}{\left|\Phi_{\text{seq,PEtab}}(\theta^*_{\text{PEtab}})\right|}$$

where $\Phi_{\text{seq,PEtab}}(\cdot)$ is the PEtab.jl-constructed objective function and $\theta^*_{\text{Exa}}, \theta^*_{\text{PEtab}}$ are the ExaModelsPetab.jl- and PEtab.jl-found optima, respectively.

We delegate the objective function evaluation of $\theta^*_{\text{Exa}}$ to $\Phi_{\text{seq,PEtab}}(\cdot)$ because the simultaneous objective $\Phi_{\text{sim}}(\cdot)$ depends explicitly on the discretized state profile, which is sensitive to collocation error. By comparing both optima from the same stiff forward solves and tolerances, we remove any bias towards false optima that may arise purely from formulation differences.

## Benchmark results

### Key observations

Across 6 of the 11 models where ExaModelsPEtab.jl achieves a solve status `0`/`0A` in Table 1, the simultaneous method on GPU is $2.3\times$ faster (geometric mean) than the sequential PEtab.jl baseline, but $6.3\times$ slower across the remaining 5/11 models ($0.68\times$ overall). The 9 failed or suboptimal models appear to occur independently of problem size, as seen in Figure 1. Comparing solvers within the simultaneous method, the GPU is faster than the CPU on all but three small models. Most notably, the simultaneous method consistently performs worse with respect to problem size exclusively on CPU but not on GPU.

**Table 1:** Benchmarking results of ExaModelsPetab.jl (GPU/CPU) and PEtab.jl. Fastest solve time is bolded, with MadNLP eligible only at status `0`/`0A`; otherwise the PEtab.jl optimizer is bolded. Bolded ROG marks a successful `0`/`0A` solve. ExaModelsPEtab.jl keys: `0:SOLVE_SUCCEEDED`, `0A:SOLVED_TO_ACCEPTABLE_LEVEL`, `0S`/`0AS` are `0`/`0A` but with ROG $> 0.02$, `T:Timeout`, `R:RESTORATION_FAILED`, `E:ERROR`. PEtab.jl keys: `0:grad_tol`, `0F:func_tol`. Best PEtab.jl optimizer time is reported: `IPN:Optim.IPNewton`, `GN:Fides.GaussNewton`, `BFGS:Fides.BFGS`.

| | ExaModelsPEtab.jl + MadNLP (GPU) | | | ExaModelsPEtab.jl + MadNLP (CPU) | | | PEtab.jl + Best optimizer (CPU) | | |
|---|---|---|---|---|---|---|---|---|---|
| Model | $t_{\text{SGM}}$ | Status | ROG | $t_{\text{SGM}}$ | Status | ROG | $t_{\text{SGM}}$ | Status | Optimizer |
| Armistead | 0.470 | 0 | **+5.1e-3** | 0.463 | 0 | **+2.3e-6** | **0.132** | 0F | GN |
| Bachmann | - | T | +7.7e-1 | - | T | +5.6e-1 | **23.336** | 0F | GN |
| Bertozzi | 0.199 | 0S | +1.9e+0 | 0.278 | 0S | +1.9e+0 | **0.189** | 0F | BFGS |
| Blasi | **0.279** | 0 | **+2.0e-7** | **0.013** | 0 | **-3.0e-7** | 0.692 | 0F | BFGS |
| Boehm | 0.329 | 0 | **-1.4e-7** | 0.548 | 0 | **-1.6e-7** | **0.009** | 0 | IPN |
| Borghans | - | R | +8.8e-1 | - | E | - | **89.965** | 0F | BFGS |
| Bruno | **0.221** | 0 | **+1.6e-8** | **0.418** | 0 | **+7.0e-11** | 0.487 | 0F | GN |
| Crauste | **0.199** | 0 | **+1.2e-3** | 0.445 | 0 | **+8.1e-6** | 0.434 | 0F | GN |
| Elowitz | - | R | +8.9e-1 | - | E | +1.0e-4 | **18.587** | 0F | GN |
| Fiedler | - | E | - | - | T | - | **0.014** | 0 | IPN |
| Laske | 21.348 | 0S | +6.5e-2 | 182.912 | 0S | +6.4e-2 | **24.660** | 0F | GN |
| Lucarelli | - | T | +2.7e-2 | - | T | +4.9e-4 | **221.146** | 0F | BFGS |
| Okuonghae | 4.724 | 0S | +5.2e-2 | 5.606 | 0AS | +1.6e-1 | **10.856** | 0F | BFGS |
| Perelson | 0.433 | 0 | **+4.7e-3** | 0.055 | 0S | +9.4e-1 | **0.135** | 0F | BFGS |
| Rahman | **0.215** | 0 | **+2.1e-7** | 0.364 | 0 | **+1.9e-7** | 0.233 | 0F | GN |
| SalazarCavazos | **39.833** | 0 | **+8.9e-4** | 481.771 | 0 | **+8.9e-4** | 215.807 | 0F | GN |
| Schwen | 29.873 | 0A | **+3.9e-6** | 115.227 | 0A | **+2.8e-7** | **6.529** | 0F | BFGS |
| Sneyd | **2.765** | 0 | **+3.4e-5** | 193.974 | 0A | **+1.2e-6** | 5.744 | 0F | BFGS |
| Zhao | 50.082 | 0AS | +8.0e-2 | 185.255 | 0AS | +1.8e-1 | **15.046** | 0F | BFGS |
| Zheng | 13.679 | 0 | **+2.9e-3** | 54.917 | 0A | **+1.0e-5** | **2.700** | 0F | BFGS |

## Trivial solves

Two of PEtab.jl's fastest solves (Boehm and Fiedler) are benchmarking artifacts resulting from the nominal guess $\theta_0$ coinciding with the optimal solution $\theta^*$. The trivial Boehm solve accounts for the smallest $0.03\times$ point in Figure 1. Here, `Optim.IPNewton` terminates on gradient stationarity after a single forward solve. Meanwhile, MadNLP GPU requires several iterations to align the forward-solved initialized state guesses with the collocation profile for primal feasibility and to converge the barrier parameters from the estimated parameter bounds. While this is a minor issue in this particular case, other stiffer models may suffer from this initialization misalignment.

## Suboptimal solutions

The 4 suboptimal solves stem from true local minima rather than extraneous problem formulation or solver errors: initializing the ExaModel at $\theta^*$ reduces the ROG to near zero. Additionally, both $\Phi_{\text{seq}}(\theta^*_{\text{Exa}})$ and $\Phi_{\text{sim}}(\theta^*_{\text{Exa}})$ evaluate to nearly the same value. While low collocation accuracy is a legitimate concern we acknowledge later, adjusting or adopting a finer mesh does not reduce the ROG for these particular models. Instead, these observations can be attributed to solution basins that differ between the simultaneous and sequential method formulations. In the future, a multi-start benchmark should be performed to delineate initial guess artifacts.

## Timeout and restoration failures

The four timeout and restoration failures fall under two modes. The two largest models Bachmann and Lucarelli (nvar $= O(10^6)$) exceed the timeout cap, stalling from primal and dual infeasibility, respectively. Lucarelli is primal feasible, yet the simultaneous objective differs greatly from the sequential objective evaluated at the terminal point, indicating a poor collocation approximation. The smaller Borghans and Elowitz models (nvar $= O(10^4)$) reach primal feasibility but terminate on restoration failure from dual infeasibility. However, a better $\theta_0$ recovers these two models, again suggesting that the simultaneous method is more sensitive to poor initialization. While PEtab.jl converges these models, they terminate on objective tolerance and also fail when forced to converge under stationarity.

## GPU kernel limitations

The Fiedler model fails to compile on GPU because the explicit analytical expressions of the steady-state profiles are embedded as initial conditions instead of utilizing the existing steady-state model build path. As a result, the callback functions generated by ExaModelsPEtab.jl preserve these massive expressions which exceed the $31.996$ KiB GPU kernel memory limit. Symbolic analysis prior to constructing the model should be considered to automatically identify proper build paths or decompose large expressions without relying on user input.

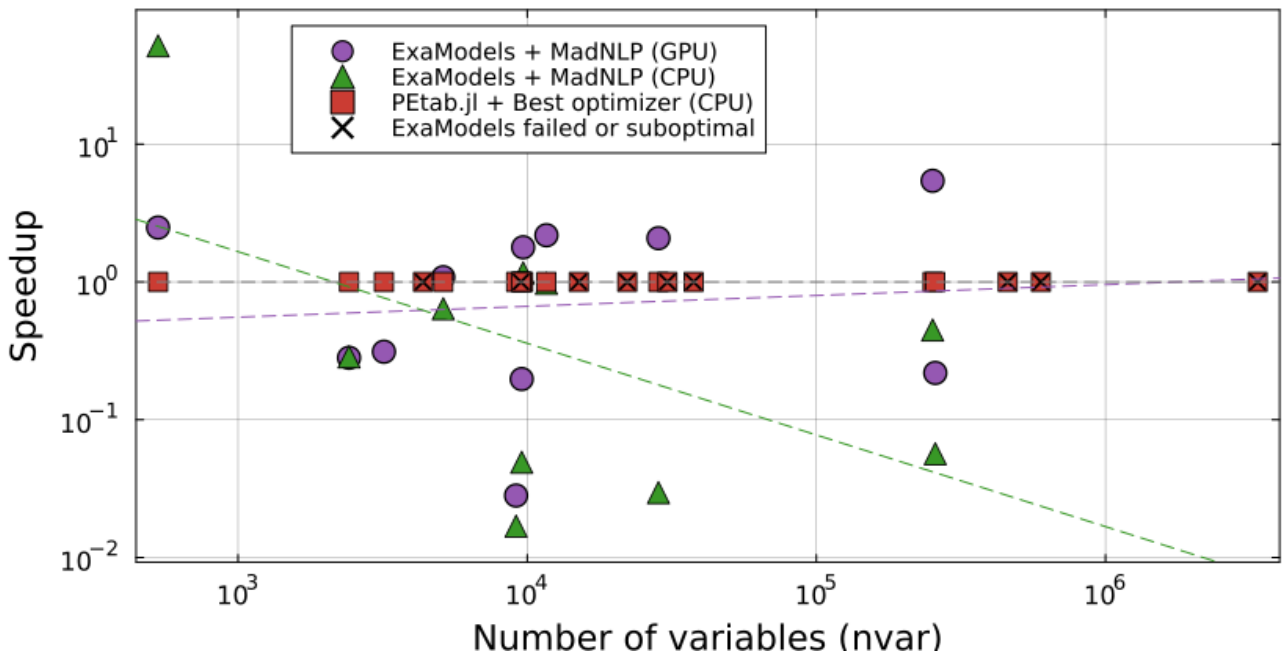

**Figure 1.** Speedup of ExaModelsPEtab.jl (GPU/CPU) vs PEtab.jl based on simultaneous method NLP size for status `0`/`0A` solves. PEtab.jl model x-axis position corresponds to its ExaModels counterpart of that size.

### Static collocation mesh limitations

Although our fixed-mesh heuristic guarantees an accurate state profile at the nominal guess $\theta_0$, this guarantee is lost as $\theta$ moves away from $\theta_0$. The initialized collocation mesh may be too coarse for accurately capturing the stiff dynamics near the optimal solution, resulting in false optima $\theta^\times$ which yield large objective errors when passed through a stiff forward solve $\Phi_{\text{sim}}(\theta^\times) \ll \Phi_{\text{seq}}(\theta^\times)$. Adaptive mesh refinement with error estimation to iteratively adjust either the interval width or $K$ may still be necessary to recover failed models that exhibit this transient stiffness. However, these methods are typically performed sequentially in an outer loop, limiting its opportunities for implementation with GPU parallelization.

## CONCLUSION

We compare a GPU-enabled simultaneous method implementation against a sequential baseline on parameter estimation benchmarks from systems biology. On models where both methods successfully solve, the simultaneous method with GPU outperforms the sequential method on roughly half of the instances. Notably, this advantage becomes more pronounced as the problem size grows only on the GPU solvers. However, this advantage is inconsistent, often failing to converge or reaching suboptimal solutions due to poor initialization, scaling issues, or discretization error from problem stiffness. Mesh refinement or better initialization strategies while maintaining full GPU compatibility should be investigated.

## ACKNOWLEDGEMENTS

This work was supported by the MIT Initiative for New Manufacturing (INM) Seed Grant.

## DECLARATION OF USE OF AI

Claude Opus 4.8 was used to debug code and assist setting up the benchmark repository.

## AUTHOR IDENTIFIERS

Author ORCIDs:
Choi JW: 0009-0000-1704-3283
Shin S: 0000-0002-9889-3278

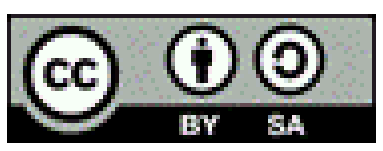